\documentclass{article}
\usepackage{mathrsfs}
\usepackage{graphicx}
\usepackage{amsmath}
\usepackage{amssymb}
 \usepackage{txfonts}

\newtheorem{thm}{Theorem}[section]

\newtheorem{lem}[thm]{Lemma}

\newtheorem{rem}[thm]{Remark}

\begin{document}

\begin{center}

{\Large \bf On $SS$-quasinormalities of the maximal subgroup series of finite
groups \footnote{* Corresponding author: J.K. Lu (jklu@gxnu.edu.cn).
W. Meng is supported by National Natural Science Foundation of China (12161021),and Guangxi Colleges and Universities Key Laboratory of Data Analysis and Computation.  J.K. Lu is supported by  Guangxi Natural Science
  Foundation Program (2024GXNSFAA010514). B. Zhang is  supported by the National Natural Science Foundation of
 China (12301022).
}
}

\end{center}

\vskip0.5cm

\begin{center}

   Wei Meng*$^a$, Jiakuan Lu$^b$, Boru Zhang$^b$

 a. School of Mathematics and Computing Science and Center for Applied Mathematics of Guangxi(GUET), Guilin University of Electronic Technology, Guilin, Guangxi, 541002, P.R. China.

 b. School of Mathematics and Statistics, Guangxi Normal University,

   Guilin, Guangxi, 541004, P.R. China,

  E-mails:   mlwhappyhappy@163.com, jklu@gxnu.edu.cn, brzhangqy@163.com

\end{center}

\begin{abstract}
Let $G$ be finite group. A subgroup $H$ of $G$ is said to be an $SS$-quasinormal
subgroup of $G$, if there exists a subgroup $B$ of $G$ such that $G = HB$ and
$H$ permutes with every Sylow subgroup of $B$.
Let $\Omega: G=G_0>G_1>\cdots>G_{n-1}>G_n=1$ be a maximal subgroup series of $G$,  where $G_i$ is a maximal subgroup of $G_{i-1}$ for every $i = 1, \ldots , n$. In this paper, we investigate the finite groups $G$ that admit an $SS$-quasinormal maximal
subgroup series, i.e.,  all $G_i$ are $SS$-quasinormal  in $G$. First, we prove that if $G$ possesses an $SS$-quasinormal maximal subgroup  series, then $G$ is solvable. Furthermore, we  show that $G$ is supersolvable if and only if $G$ possesses an $SS$-quasinormal maximal subgroup  series which is  subnormal in $G$.
\end{abstract}




\section{Introduction}
Throughout this paper, all groups are finite. The investigation of how subgroup properties affect the structure of finite groups is a prominent topic in finite group theory, and this line of research has yielded numerous results.  By taking $X$ to be a subset of some subgroups of $G$ and requiring that every subgroup in $X$ has "well-behaved" properties, our goal is to investigate the structure of group $G$.  Specifically, in  existing  research, $X$ is frequently taken to the set of maximal subgroups of $G$, the set of minimal subgroups of $G$, or the set of subgroups of $G$ with a specific order (see \cite {SK}). Moreover, from the viewpoint of $G$'s subgroup lattice, such subsets $X$ usually constitute a horizontal cross-section between $G$ and the identity subgroup $1$.  In line with this perspective, let $G$ be a group and
\begin{center}
$\Omega: G=G_0>G_1>\cdots>G_{n-1}>G_n=1$
\end{center}
be a maximal subgroup series of $G$, meaning that $G_i$ is a maximal subgroup of $G_{i-1}$ for every $i = 1, \ldots , n$.
 Some interesting results can be obtained under assuming that the maximal subgroup series has certain normality. For example, the series $\Omega$ is said to be central in $G$ if $[G,G_{i-1}] \leq  G_i$ for every $i = 1, \ldots, n$; and it is said to be
normal (or subnormal) in $G$ if all $G_i$ are normal (or subnormal) in $G$. The following results are well-known:

(1)  $G$ is nilpotent if and only if $G$ possesses a maximal subgroup series that is central in $G$.

(2) $G$ is supersolvable if and only if $G$ possesses a maximal subgroup series that is normal in $G$.

(2) $G$ is solvable if and only if $G$ possesses a maximal subgroup series that is subnormal in $G$.

Recall that a subgroup $A$ of $G$ is said to be $S$-permutable (or $S$-quasinormal)  in $G$ if $A$ is permutable with all Sylow
subgroups of $G$(see \cite {KE}).  A subgroup $A$ of $G$ is called completely conditional permutable (abbreviated as $c$-$c$-permutable)
in $G$ if, for each subgroup $H$ of $G$, there exists an element $g\in\langle A,H\rangle$ such that $AH^g = H^gA$ (see\cite {GU,QI1}).

Recently, Qian and Tang \cite {QI2} studied the finite groups $G$ that admit an $S$-permutable (or $c$-$c$-permutable) maximal
subgroup series, i.e., a series $G = G_0 > \cdots > G_i > \cdots> G_n = 1$ where all $G_i$ are $S$-permutable  (or $c$-$c$-pemutable) in $G$.
 They proved the following theorem:
 \begin{thm} \cite [Theorems 1,2] {QI2}
 $G$ is supersolvable if and only if $G$ possesses an $S$-permutable (or $c$-$c$-pemutable) maximal subgroup series.
 \end{thm}

In 2008,  Li \cite {LI1} generalized $S$-permutable   subgroups to
$SS$-quasinormal subgroups. A subgroup $H$ of $G$ is said to be an $SS$-quasinormal
subgroup of $G$, if there exists a subgroup $B$ of $G$ such that $G = HB$ and
$H$ permutes with every Sylow subgroup of $B$.   Some interesting
results of $SS$-quasinormal subgroups can be found  in \cite {LI1,LI2}.

In the light of above investigations, we concern  the finite groups $G$ that admit an $SS$-quasinormal maximal
subgroup series, i.e., a series $G = G_0 > \cdots> G_i >\cdots> G_n = 1$ where all $G_i$ are $SS$-quasinormal  in $G$.
Our first goal is to prove the following theorem.
\begin{thm}
If $G$ possesses an $SS$-quasinormal maximal subgroup  series, then $G$ is solvable.
\end{thm}

The following remark shows that if $G$ only possesses an $SS$-quasinormal maximal subgroup series, then $G$ is not necessarily a supersolvable group.
\begin{rem}
Let $G=A_4$ and $H$ be a subgroup of $G$ of order $3$. Then
\begin{center}
 $\Omega_1: G>H>1$
  \end{center}
  is an $SS$-quasinormal maximal subgroup series of $G$. However, $G$ is not supersolvable. On the other hand, let $K$ be a subgroup of $G$ of order $4$ and $L$ be any subgroup of $K$ of order $2$. Then
  \begin{center}
  $\Omega_2:G>K>L>1$
  \end{center}
   is an subnormal maximal subgroup series of $G$. But $\Omega_2$ is not $SS$-quasinormal in $G$.
\end{rem}

Inspired by Remark 1.3, our second goal is to show   the following result.
\begin{thm}
$G$ is supersolvable if and only if $G$ possesses an $SS$-quasinormal maximal subgroup series which is  subnormal in $G$.
\end{thm}

All unexplained notations and terminologies
    are standard and can be found in \cite {DO, HB}.

\section{Preliminaries}

In this section, we collect some results which will be used in
the proof of the main results.

\begin{lem} \cite [Lemma 2.1] {LI1}
Suppose that $H$ is an $SS$-quasinormal subgroup of $G$, $K\leq G$ and $N$ is a normal subgroup of $G$. Then, we have the
following:

(1) If $H \leq K$, then $H$ is an $SS$-quasinormal subgroup of $K$.

(2) $HN/N$ is an $SS$-quasinormal subgroup of $G/N$.
\end{lem}

\begin{lem} \cite [Lemma 2.2] {LI1}
  Let $H$ be a nilpotent subgroup of $G$.
Then, the following statements are equivalent:

(1) $H$ is an $S$-quasinormal subgroup of $G$.

(2) $H \leq F(G)$ and $H$ is an $SS$-quasinormal subgroup of $G$.
\end{lem}

\begin{lem}
Let $M$ be a maximal subgroup of $G$. If $M$ is $SS$-quasinormal in $G$, then $|G:M|$ is a prime power.
\end{lem}
{\noindent\bf  Proof.} Suppose that there are two distinct primes, says, $p$ and $q$ such that $pq\mid |G:M|$. Since $M$ is $SS$-quasinormal in $G$, there there exists a subgroup $B$ of $G$ such that $G = MB$ and
$H$ permutes with every Sylow subgroup of $B$, in particular, $MB_p=B_pM$ is a proper subgroup of $G$, where $B_p$ is a Sylow $p$-subgroup of $B$. So we have   $M<MB_p<G$. It  contradicts that $M$ is a maximal subgroup of $G$.   So $|G:M|$ is a prime power. $\Box$

\begin{lem}
Let $N$ be a normal subgroup of $G$ and assume that $G$ possesses an $SS$-quasinormal   maximal subgroup series. Then $G/N$
also has an $SS$-quasinormal  maximal subgroup series.
\end{lem}
{\noindent\bf  Proof.} Let $G = G_0 > \cdots> G_i > \cdots > G_n = 1$ be an $SS$-quasinormal  maximal subgroup series of $G$.
Since $G_i$ is $SS$-quasinormal   in $G$,  w $G_iN/N$ is also $SS$-quasinormal  in $G/N$ by Lemma 2.1 (2). Write
$\overline {G} = G/N $ and $\overline {G_i} = G_iN/N$. Let us investigate the following subgroup series of $\overline {G}$:
\begin{center}
$\overline{G} = \overline {G_0} \geq\cdots\geq\overline{G_i}\geq\cdots\geq\overline{G_n} = 1$.
\end{center}
For each $i = 1,\ldots , n$, we see that either $\overline{G_i} = \overline{G_{i-1}}$ or $\overline{G_i}$ is maximal in $\overline{G_{i-1}}$. Therefore, after removing
the equal terms in the above series, we obtain an $SS$-quasinormal   maximal subgroup series of $\overline{G}$. $\Box$

\begin{lem} \cite [Theorem 1] {GUR}
Let $G$ be a non-abelian simple group. If $H$ is a proper subgroup of $G$ with index $p^a$, where $p$ is a prime, then one of the following holds:

(1) $G= A_n$ and $H=A_{n-1}$, $n=p^a$;

(2) $G=PSL(n,q)$ and  $H$ is the stabilizer of a line or hyperplane, $|G:H|=(q^n-1)/(q-1)=p^a$ and $n$ is a prime;

(3) $G=PSL(2,11)$, $H=A_5$;

(4) $G=M_{23}$ and $H=M_{22}$;

(5) $G=M_{11}$ and $H=M_{10}$;

(6) $G=PSU(4,2)$ and $|G:H|=27$.
\end{lem}

\begin{lem} \cite [Lemma 2] {BE}
Let $ N = S_1\times\cdots\times S_t$ be a direct product of isomorphic non-abelian simple groups, and let
$M$ be a maximal subgroup of $N$ with $S_1\not\leq M$.  Then one of the following assertions holds:

(1) $M = D\times S_2\times\cdots\times S_t$, where $D$ is maximal in $S_1$.

(2) One of the subgroups $S_2,\ldots, S_t$, say $S_2$, is not contained in $M$, then $M = D \times  S_3\times\cdots\times S_t$, where
$D \cap S_1 = D \cap S_2 = 1$ and $S_1\cong S_2\cong D < S_1\times S_2$.
\end{lem}
\section{Proof  of Theorem}

{\noindent\bf  Proof of Theorem 1.2.}  Assume that $G$ is non-solvable and possesses  an $SS$-quasinormal maximal subgroup series:
\begin{center}
$\Omega: G=G_0>G_1>\cdots>G_{n-1}>\cdots>G_n=1$.
\end{center}
Let $N$ be a minimal normal subgroup of $G$. Then $G/N$ also has an $SS$-quasinormal
maximal subgroup series by Lemma 2.4. It follows by induction that $G/N$ is  solvable. If $G$ has distinct minimal
normal subgroups, say $N_1$ and $N_2$, then both $G/N_1$ and $G/N_2$ are  solvable, and so is $G$.
Consequently, we may assume that  $G$ possesses a unique minimal normal subgroup, say
$N$. Since $G$ is non-solvable, $N$ is a direct
product of some isomorphic non-abelian simple groups  and $C_G(N)=1$.

Applying Lemma 2.1(1), we know that  $G_1>\cdots>G_{n-1}>\cdots>G_n=1$ is an $SS$-quasinormal maximal subgroup series of $G_1$. It follows by induction that $G_1$ is solvable.  So  $N\not\leq G_1$ and hence $G_1\cap N$ is a proper subgroup of $N$.  Observe that $|G_{i-1}:G_i|$ is a prime power by Lemma 2.3,  and that $|N \cap G_{i-1} : N\cap G_i|$ divides $|G_{i-1} : G_i|$. We have  $|N \cap G_{i-1} : N\cap G_i|$ is a prime power or $1$ for every $i = 1,\ldots , n$. In particular, $|N:N\cap G_1|=|N\cap G_0:N\cap G_1|=p^{\alpha}$ for some prime $p$.  This shows that $N$ possesses a  solvable   subgroup with index prime power, namely $N\cap G_1$. By Lemma 2.6, we get that $N$ is a non-abelian simple group.  Moreover, $N$ will be isomorphic to one of the group in Lemma 2.5.

Let $j$ be the largest index such that $|N : N \cap G_j|$ is a $p$-power. Let $B = N\cap G_j$ and $C = N \cap G_{j+1}$. Then $|B:C|=q^{\beta}$,  where $q$ is a prime different from $p$. Since $G_{j+1}$ is $SS$-quasinormal in $G$, there exists a subgroup $B$ of $G$ such that $G=G_{j+1}B$ and $G_{j+1}B_p=B_pG_{j+1}$, where $B_p\in Syl_p(B)$. This implies that $G_{j+1}B_p$ is a proper subgroup of $G$ and $|G|_p=|G_{j+1}B_p|_p$.  Consequently, there exists a subgroup $G_p\in Syl_p(G)$ such that $G_{j+1}G_p=G_pG_{j+1}$.  It follows that
\begin{align*}
 C(N\cap G_p)&=(N\cap G_{j+1})(N\cap G_p)\\
 &=N\cap(G_{j+1}(N\cap G_p))\\
 &=  (G_{j+1}(N\cap G_p))\cap N\\
 &=((N\cap G_p)G_{j+1})\cap N)\\
 &=(N\cap G_p)C.
\end{align*}
Set $N_p=N\cap G_p$, then $N_p\in Syl_p(N)$. So
\begin{center}
$|N:CN_p|=|N:CN_P|_{p'}=|N:C|_{p'}=|B:C|=q^{\beta}$.
\end{center}

Now, the non-abelian simple group $N$ admits subgroups $G_1 \cap N$ and $CN_p$  such that $|N : G_1 \cap N|$ and
$|N : CN_p|$ are distinct prime powers. Note that   PSL$(2,7)$ is the only simple group with
subgroups of two different prime power indices (see \cite {GUR}), so $N\cong$PSL$(2,7)$.

Finally, since $N$ is the unique minimal normal subgroup $G$ and $C_G(N)=1$, by N/C-theorem,
\begin{center}
$N\leq G=G/C_G(N)\precapprox$Aut$(N)$.
\end{center}
Consequently, since Aut(PSL$(2,7)) \cong$PSL$(2,7).Z_2$, we have
\begin{center}
PSL$(2,7)\cong N\leq G\precapprox$ Aut(PSL$(2,7))\cong$ PSL$(2,7).Z_2$.
\end{center}
So $G=N\cong$PSL$(2,7)$ or $G\cong$PSL$(2,7).Z_2$.
By checking the maximal  subgroup series  of $G$, we can  conclude that $G$ doesn't possesses $SS$-quasinormal maximal series.
The proof of theorem is complete. $\Box$\\

{\noindent\bf  Proof of Theorem 1.4.} Suppose that $G$ is supersolvable. Then every chief series of $G$ is necessarily an
$SS$-quasinormal maximal subgroup series of $G$ and is subnormal in $G$.  Conversely, suppose that the theorem is not true
and let $G$ be a counterexample of the smallest order. Let
\begin{center}
$\Omega: G=G_0>G_1>\cdots>G_{n-1}>G_n=1$.
\end{center}
be a subnormal maximal subgroup series of $G$,  and assume that $\Omega$ is $SS$-quasinormal in $G$.
 It is clearly that $G$ is solvable and  $G_i$ is normal in $G_{i-1}$ with prime index. By Lemma 2.1(1), each $G_i$   satisfies the hypothesis of theorem, so $G_i$ is supersolvable for $i\geq 1$ by choice of $G$.

Let $N$ be a minimal normal subgroup of $G$. Then  $G/N$ satisfies the hypothesis of theorem by  Lemma 2.4 and hence  $G/N$   is supersolvable. Consequently, if $G $ has two distinct minimal normal subgroups, says $N_1$ and $N_2$, then both $G/N_1$ and $G/N_2$ are supersolvable, and so is $G$. This contradicts to the choice of $G$.  Therefore, $G$ possesses a unique minimal normal subgroup, says $N$. Since $G$ is solvable, we may assume that $N$ is an elementary abelian $p$-group of order $p^e$ for some prime $p$.
Furthermore, suppose that $\Phi(G) > 1$, then $N\leq\Phi(G)$ and hence $G/\Phi(G)$ is supersolvable.  Thus, $G$ is supersolvable,  which is another contradiction.  Therefore,
 we may assume that $\Phi(G) = 1$.
  Moreover, applying the solvability of $G$ again, there is  a maximal subgroup $H$ of $G$ such that $G = HN = H\ltimes N$,
where $H\cong G/N$  is supersolvable.  Now, it is easy to see that $N = O_p(G) = F(G) = C_G(N)$, and $C_H(N) = 1$.

Furthermore, observe that $G_1$ is normal in $G$ and $G_{n-1}$ is a subnormal subgroup of $G$ of prime order.  We have
\begin{center}
$G_{n-1}\leq N=O_p(G)=F(G)<G_1$.
\end{center}
Note that $G_{n-1}$ is $SS$-quasinormal in $G$. Applying Lemma 2.2, $G_{n-1}$ is $S$-permutable in $G$.

On the other hand, we have  $G_1=G_1\cap G=G_1\cap HN=(G_1\cap H)\ltimes N$. Set $H_1=G_1\cap H$, then $H_1$ is normal in $H$ as $G_1$ is normal in $G$.
  By the supersolvability of $G_1$, we have $O_{p',p,p'}(G_1)=G_1$. Furthermore, since $O_{p'}(G_1)\leq O_{p'} (G) = 1$, $G_1$
has a normal Sylow $p$-subgroup which is also normal in $G$. Consequently, $N $is exactly the normal Sylow
$p$-subgroup of $G_1$. Hence, $H_1$ is a $p'$-group.

We claim that $H$ is also a $p'$-group. If not,  then $H = P_0\ltimes H_1$, where $P_0 \in Syl_p(H)$ has order $p$.  By the
supersolvability of $G_1$, we conclude that $H_1\cong G_1/N=G_1/C_{G_1}(N)$
  is an abelian group with exponent
dividing $p-1$. In particular, $p$ is the largest prime divisor of $|H|$. Since $H$ is  also supersolvable, $H$ has a normal
Sylow $p$-subgroup. This leads to $O_p(G) = P_0N > N$, a contradiction. Hence, $H$ is a $p'$-group, as desired.

 Finally, observe that $G_{n-1}$ is $S$-permutable in $G$,   It follows that $G_{n-1}H=HG_{n-1}\leq G$. Since $H$ is a maximal subgroup of $G$, we have $G=HG_{n-1}$ which implies $|G_{n-1}|=|G:H|=|N|$. Therefore,  $G$
supersolvable. This is a final contradiction. The proof of theorem is complete.$\Box$



\bibliographystyle{amsplain}

\end{document}